\documentclass[10pt,twocolumn,twoside]{IEEEtran}
\usepackage{cite}
\usepackage{amsmath,amssymb,amsfonts}
\usepackage{graphicx}
\usepackage{textcomp}
\usepackage{xcolor}
\usepackage{amsmath}
\usepackage{algorithm}
\usepackage{stfloats}
\usepackage{blindtext}
\usepackage{multicol,lipsum}
\usepackage{algpseudocode}
\usepackage{lipsum}
\usepackage{ulem}
\ifCLASSINFOpdf
\else
\fi
\begin{document}
	%
	\title{Reciprocity of Algorithms Solving Distributed Consensus-Based Optimization and Distributed Resource Allocation}
	%
	%
	%
	
	\author{Seyyed~Shaho~Alaviani,~\IEEEmembership{Member,} Atul~Gajanan~Kelkar,~\IEEEmembership{Senior Member,}
		and~Umesh~Vaidya,~\IEEEmembership{Senior Member,~IEEE}
		\thanks{Seyyed Shaho Alaviani is with the Mechanical \& Industrial Engineering Department, University of Minnesota, Duluth, MN 55812 USA e-mail: salavian@umn.edu.}
		\thanks{Atul Gajanan Kelkar and Umesh Vaidya are with the Department of Mechanical Engineering, Clemson University, Clemson, SC, 29634 USA e-mails: atul@clemson.edu, uvaidya@clemson.edu.}
		\thanks{A preliminary version of this paper has appeared without proofs in \cite{alavianiMED2021}. This work has been done while Seyyed Shaho Alaviani was at Clemson University, Clemson, South Carolina, USA. }
	}
	
	%
	%

	\markboth{}%
	{Alaviani \MakeLowercase{\textit{et al.}}: Recipr. Alg. Solv. Dist. Consen. Optim. \& Dist. Res. Alloc.}
	%



	\maketitle
	
	\begin{abstract}
		This paper aims at proposing a \textit{procedure} to derive distributed algorithms for distributed consensus-based optimization by using distributed algorithms for network resource allocation and vice versa over \textit{switching} networks with/without synchronous protocol. It is shown that \textit{first-order} gradient distributed consensus-based optimization algorithms can be used for finding an optimal solution of distributed resource allocation with synchronous protocol under \textit{weaker} assumptions than those given in the literature for non-switching (static) networks. It is shown that \textit{first-order} gradient distributed resource allocation algorithms can be utilized for finding an optimal solution of distributed consensus-based optimization. The results presented here can be applied to \textit{time-varying} and \textit{random directed} networks \textit{with} or \textit{without} synchronous protocol with \textit{arbitrary} initialization. As a result, several algorithms can now be used to derive distributed algorithms for \textit{both} consensus-based optimization and resource allocation, that can overcome limitations of the existing results. While the focus of this paper is on the first-order gradient algorithms, it is to be noted that the results also work with \textit{second-order} gradient algorithms.
	\end{abstract}
	
	\begin{IEEEkeywords}
		distributed optimization, resource allocation, random networks, asynchronous. 
	\end{IEEEkeywords}

	%
	\IEEEpeerreviewmaketitle

	\section{Introduction}
	Resource Allocation problem appears in a wide spectrum of domains including higher education \cite{highereducation}, economics \cite{economics}, and health care \cite{healthcare}, to name a few. \textit{Resource Allocation} problem is to allocate the network resources 
	among a group of agents while optimizing a certain performance
	index. Network resource allocation is also fundamental and an important
	problem that arises in a variety of application domains
	within engineering such as the media access control in communication networks \cite{communi}, signal processing \cite{signalpro}, the load demand management \cite{loaddemand}, and economic dispatch problems in power systems \cite{economicdisp} (see the survey paper \cite{survey} for more details). 
	
	In \textit{centralized} allocation schemes \cite{centralized}, the allocation decisions are made centrally by gathering all network data together and coming up with a decision to satisfy certain central objective and then sending the decision back to the agents. However, centralized algorithms may \textit{not} be effective and practical in a large-scale network due to heavy communication, computational burden, private information, or complicated network structure. Therefore, designing \textit{decentralized} and/or \textit{distributed} algorithms for resource allocation is very important and highly desirable. Several researchers have focused their attention to develop decentralized and/or distributed algorithms \cite{economicdisp}, \cite{dist1}-\cite{distrandom2}. Among these, the authors of \cite{distbconnectivity1}-\cite{distrandom2} have considered \textit{switching} graphs and imposed assumptions such as connectivity at any time, B-connectivity\footnote{There exists a bounded time interval such that the union of the graphs is strongly connected, and each edge transmits a message at least once.}, independent and identically distributed (i.i.d.), or Markov chain on the underlying graph. In \cite{PEVresourceallo3} and \cite{distbconnectivity3}, \textit{asynchronous}\footnote{In a \textit{synchronous} protocol, all nodes activate at the same time and perform communication updates. On the other hand, in an \textit{asynchronous} protocol, each node has its concept of time defined by a local timer, which randomly triggers either by the local timer or by a message from neighboring nodes. The algorithms guaranteed to work with no \textit{a priori} bound on the time for updates are called \textit{totally asynchronous}, and those that need the knowledge of \textit{a priori} bound, known as B-connectivity assumption, are called \textit{partially asynchronous} (see \cite{17} and \cite[Ch. 6-7]{qlearnbertsekas}). As the dimension of the network increases, synchronization becomes an issue.  } algorithms have been proposed for solving distributed resource allocation problem. {\textit{Hence, only few results exist for decentralized and/or distributed resource allocation algorithms over switching graphs or with asynchronous protocol, where the proposed algorithms are first-order gradient.}   The authors of \cite{secondorder1}-\cite{secondorder2} have proposed discrete-time distributed \textit{second-order} algorithms over static (non-switching) graphs under synchronous protocol for resource allocation optimization. The second-order gradient algorithms provide faster convergence rates.

	Another challenge in networked systems is distributed \textit{consensus-based} optimization. Distributed consensus-based optimization appears in many problems such as sensor networks or data regression (see survey papers  \cite{survey2019}-\cite{nedic} for more applications). The aim of agents in distributed consensus-based optimization is to get an \textit{agreement} on some value that is an optimal solution to the problem. Note that the agents may not \textit{necessarily} get an agreement on some value in distributed resource allocation. As seen in surveys \cite{survey2019}-\cite{nedic}, several distributed algorithms for consensus-based optimization have been proposed over \textit{static}, \textit{time-varying}, or \textit{random} networks \textit{with} or \textit{without} asynchronous protocols. Recently, the \textit{first} totally asynchronous (see footnote 2) first-order gradient algorithm has been proposed in \cite{alavianiTAC} for distributed consensus-based optimization. \textit{To the best of authors' knowledge, a totally asynchronous algorithm has not been proposed for distributed resource allocation. }

	
	Upon the review of the literature on \textit{distributed} consensus-based optimization and \textit{distributed} resource allocation and their importance in engineering, the natural questions that arise are: \textit{Question 1)} is there a relation between these two problems?; \textit{Question 2)} is it possible to use algorithms for distributed consensus-based optimization for solving distributed resource allocation problem?; and \textit{Question 3)} is it possible to utilize algorithms for distributed resource allocation for solving distributed consensus-based optimization problem? Although it has been established that there is a relation between these two problems in one direction (reviewed below), it is not known if there exist two-way reciprocity for solving these problems. \textit{This paper is addressing this question. }

	

	With regard to \textit{Question 1}, it is known that the \textit{Lagrange's dual} problem of resource allocation is consensus-based optimization, e.g., see \cite[Sec. III.A]{Zhangconsensusresource} for proof; and \textit{Fenchel's dual} problem \cite{dualitybook} of distributed consensus-based optimization is resource allocation, e.g., see \cite[Sec. III.A]{Fencheldualityproof} for proof. It needs to be mentioned that Fenchel's dual problem is a \textit{duality} problem (like Lagrange's duality) represented by (Fenchel) conjugate functions (e.g., see \cite{dualitybook} for details). Note that solving both Fenchel's and Lagrange's dual problems is a \textit{challenging} problem since the cost function in dual problems is \textit{not} an \textit{explicit} function of dual variables and, in fact, is determined by another minimization problem (see Comment 1 is Section IV for details).
	
	The existing literature presents some answers to \textit{Question 2} which are explained next. In \cite{changconsresour}, it is mentioned that a distributed consensus-based optimization algorithm may be used for resource allocation through solving the dual of resource allocation; however, the authors of \cite{changconsresour} did not present any method to do so. Later, the authors in \cite{Nedichconsensusresource} showed the relation that the \textit{first-order} optimality conditions of the two aforementioned problems have when underlying graph is \textit{static}, connected, and \textit{undirected} with \textit{synchronous} protocol and when the cost function of each agent is proper, closed, and convex. This leads to the use of decentralized or distributed consensus-based optimization \textit{first-order} gradient algorithms with arbitrary initialization for solving resource allocation \textit{without} investigating the Lagrangian dual relationship between the two problems. Note that the optimality conditions given in \cite{Nedichconsensusresource} \textit{only} work for static networks. Recently, it was shown in \cite{Zhangconsensusresource} how to use distributed consensus-based \textit{first-order} gradient optimization algorithms for solving resource allocation over \textit{static} \textit{directed} networks with \textit{synchronous} protocol when cost function of each agent is \textit{strongly} convex. Based on the explanation above, one can apply existing distributed consensus-based optimization algorithms \textit{only} over \textit{static} networks with synchronous protocols, and \textit{we cannot derive algorithms for distributed resource allocation over switching networks and/or with asynchronous protocol to overcome restrictions of existing algorithms. }  
	

	\textbf{Contribution:} In this paper, we present how to utilize distributed consensus-based optimization algorithms to find an optimal solution of a distributed resource allocation problem and vice versa over \textit{switching} networks with/without synchronous protocol. Therefore, the objective of this paper is \textit{not} to propose a distributed algorithm and to analyze its convergence; in fact, we propose a \textit{procedure} (i.e., \textbf{P1}-\textbf{P3}) to derive a distributed algorithm for consensus-based optimization and resource allocation while its convergence analysis has already been analyzed by imposing suitable assumptions. As such, several algorithms can be employed, by imposing suitable assumptions, to derive distributed algorithms for \textit{both} consensus-based optimization and resource allocation, that can overcome limitations of the existing results. We show that \textit{first-order} gradient algorithms for  distributed consensus-based optimization  can also be used to find an optimal solution of distributed resource allocation over static directed networks with synchronous protocol under \textit{weaker} assumptions (than strong convexity in \cite{Zhangconsensusresource})  such as strict convexity, Lipschitz gradient, or coercivity. We show that distributed resource allocation \textit{first-order} gradient algorithms can also be utilized to find an optimal solution of distributed consensus-based optimization. Thus, our results show the applicability of distributed resource allocation algorithms for solving distributed consensus-based optimization, namely to present an affirmative answer to \textit{Question 3}. It is important to emphasize that \textit{conditions} given here can lead to use distributed algorithms for both consensus-based optimization and resource allocation over \textit{static}, \textit{time-varying}, or \textit{random} \textit{directed} networks \textit{with} or \textit{without} asynchronous protocols under \textit{arbitrary} initialization. While the technical issue behind is challenging (see the fifth paragraph in this section) to demonstrate how significant the improvement is made when fixing the gap in the literature, we apply the first totally asynchronous (see footnote 2) first-order gradient algorithm \cite{alavianiTAC} intended for distributed consensus-based optimization to derive the \textit{first} totally asynchronous first-order gradient-based algorithm for distributed resource allocation (which is a missing gap in the literature) where the algorithm does \textit{not} require \textit{a priori} B-connectivity (see footnote 1) or distribution assumption on switching communication graphs for convergence. The derivation of first-order distributed consensus-based optimization algorithms from distributed resource allocation algorithms is similar to the derivation of distributed resource allocation algorithms from distributed consensus-based optimization algorithms, namely the \textbf{P1}-\textbf{P3} are still valid. While we focus in this paper on \textit{first-order} gradient algorithms, \textit{the results also work with \textit{second-order} gradient (see Lemma 1 and Observation 1) algorithms where the related results are given in \cite{alavianisecondoreder}.} \textit{This paper provides proofs and related details for theorem given in [1] as they were omitted due to space limitations.}

	This paper is organized as follows: In Section II, necessary preliminaries on operators and stochastic convergence are provided. Formulations of the resource allocation and consensus-based optimization problem are given in Section III. In Section IV, the procedure for using distributed consensus-based optimization algorithms to solve distributed resource allocation is given. Section V presents how to use distributed resource allocation's algorithms to solve distributed consensus-based optimization.

	\textit{Notations:} $\Re$ denotes the set of all real  numbers. $(.)^{T}$ represents the transpose of a vector or a matrix. For any vector $z \in \Re^{n}, \Vert z \Vert_{2}=\sqrt{z^{T}z},$ and for any matrix $Z \in \Re^{n \times n}, \Vert Z \Vert_{2}=\sqrt{\lambda_{max}(Z^{T}Z)}=\sigma_{max}(Z)$ where $\lambda_{max}$ represents the maximum eigenvalue, and $\sigma_{max}$ represents the largest singular value. Sorted in an increasing order with respect to real parts, $\lambda_{2}(Z)$ represents the second eigenvalue of a matrix $Z$. $Re(r)$ represents the real part of the complex number $r$. $I_{n}$ represents Identity matrix of size $n \times n$ for some $n \in \mathbb{N}$ where $\mathbb{N}$ denotes the set of all natural numbers. $\otimes$ represents Kronecker product. $E[x]$ denotes Expectation of the random variable $x$. $\nabla f(x)$ represents the gradient of the function $f$ at $x$. $\partial f(x)$ represents the sub-differential of the function $f$ at $x$. $\textbf{0}_{n}$ and $\textbf{1}_{n}$ denote the vector of dimension $n$ whose elements are all zero and one, respectively. 
	

	\section{Preliminaries}
	A vector $v \in \Re^{n}, n \in \mathbb{N},$ is said to be a \textit{stochastic vector} when its components $v_{i}, i=1,2,...,n$, are non-negative and their sum is equal to 1; a square $n \times n$ matrix $V$ is said to be a \textit{stochastic matrix} when each row of $V$ is a stochastic vector. A square $n \times n$ matrix $V$ is said to be \textit{doubly stochastic} when both $V$ and $V^{T}$ are stochastic matrices.

	Let $\mathcal{H}$ be a real Hilbert space with norm $\Vert . \Vert $ and inner product $\langle .,. \rangle$. An operator $A:\mathcal{H} \longrightarrow \mathcal{H}$ is said to be \textit{monotone} if $\langle x-y,Ax-Ay \rangle  \geq 0$ for all $x,y \in \mathcal{H}$. $A:\mathcal{H} \longrightarrow \mathcal{H}$ is called $\rho$\textit{-strongly monotone} if $\langle x-y,Ax-Ay \rangle \geq \rho \Vert x-y \Vert^{2}$ for all $x,y \in \mathcal{H}$. A function $f(.)$ is $\rho$\textit{-strongly convex} if $\langle x-y,\nabla f(x)-\nabla f(y) \rangle \geq \rho \Vert x-y \Vert^{2}$ for all $x,y \in \mathcal{H}$. Therefore, a function is $\rho$-strongly convex if its gradient is $\rho$-strongly monotone. A function $f(.)$ is \textit{strictly convex} if $\forall x,y, x \neq y, \forall \kappa \in (0,1)$, we have 
	$$f(\kappa x +(1-\kappa) y) < \kappa f(x)+(1-\kappa) f(y).$$
	
	A mapping $B:\mathcal{H} \longrightarrow \mathcal{H}$ is said to be $K$\textit{-Lipschitz continuous} if there exists a $K > 0$ such that 
	$\Vert Bx-By \Vert \leq K \Vert x-y \Vert$ for all $x,y \in \mathcal{H}$. Let $S$ be a nonempty subset of a Hilbert space $\mathcal{H}$ and $Q:S \longrightarrow \mathcal{H}$. The point $x$ is called a \textit{fixed point} of $Q$ if $x=Q(x)$.

	Let $(\Omega^{*},\sigma)$ be a measurable space ($\sigma$-sigma algebra) and $C$ be a nonempty subset of a Hilbert space $\mathcal{H}$. A mapping $x:\Omega^{*} \longrightarrow \mathcal{H}$ is \textit{measurable} if $x^{-1}(U) \in \sigma$ for each open subset $U$ of $\mathcal{H}$. The mapping $T:\Omega^{*} \times C \longrightarrow \mathcal{H}$ is a \textit{random map} if for each fixed $z \in C$, the mapping $T(.,z):\Omega^{*} \longrightarrow \mathcal{H}$ is measurable, and it is \textit{continuous} if for each $\omega^{*} \in \Omega^{*}$ the mapping $T(\omega^{*},.):C \longrightarrow \mathcal{H}$ is continuous.

	Let $(\Omega,\mathcal{F},\mathbb{P})$ be a complete probability space where $\Omega$ denotes the sample space, $\mathcal{F}$ denotes a $\sigma$-algebra on $\Omega$, and $\mathbb{P}$ is a probability measure on $\mathcal{F}$. A sequence of random variables $x_{t}$ is said to
	\begin{itemize}
		\item \textit{converge almost surely} to $x$ if there exists a subset $A \subseteq \Omega$ such that $\mathbb{P}(A)=0$, and for every $\omega \notin A$, $\lim_{t \longrightarrow \infty} \Vert x_{t}(\omega)-x(\omega) \Vert=0.$
		
		\item \textit{converge in mean square} to $x$ if $E [\Vert x_{t}-x \Vert^{2}] \longrightarrow 0$ as $t \longrightarrow \infty.$
	\end{itemize}

	\textbf{Definition 1} \cite{alavianiTAC}: If there exists a point $\hat{x} \in \mathcal{H}$ such that $\hat{x}=T(\omega^{*},\hat{x})$ for all $\omega^{*} \in \Omega^{*}$, it is called \textit{fixed-value point}, and $FVP(T)$ represents the set of all fixed-value points of $T$.

	\textbf{Definition 2} \cite{dualitybook},\cite{boydbook}: The \textit{(Fenchel) conjugate} of a function $f: \Re^{n} \longrightarrow \Re$ is defined as 
	$$f^{*}(y) := \sup_{x}  (y^{T} x- f(x)).$$
	Moreover, the conjugate function is a convex function.
	
	\textbf{Definition 3} \cite{dualitybook}: A function $f: \Re^{n} \longrightarrow \Re$ is said to be \textit{$0-$coercive} and \textit{$1-$coercive} when $$\displaystyle \lim_{\Vert x \Vert \longrightarrow +\infty} f(x)=+\infty \quad{}\text{and} \quad{} \displaystyle \lim_{\Vert x \Vert \longrightarrow +\infty} \frac{f(x)}{\Vert x \Vert}=+\infty,$$ respectively.



	\section{Problem Formulation}

	We use the \textit{switched}\footnote{Switched systems can be roughly
		divided into two groups: those subject to \textit{arbitrary}, i.e., state-independent switching, and those subject to \textit{state-dependent} switching (see \cite{11222} for details). We require to mention that we consider \textit{arbitrary} switching in this paper. Note that \textit{weighted matrix of the graph, Laplacian,} or \textit{Adjacency} matrix can be used for information of the network in our formulation. For instance, we use weighted matrix of the graph in this paper. Laplacian has been utilized in \cite{alavianisecondoreder} for second-order gradient algorithms.} network topology similar to the one used in \cite{alavianiTAC} wherein a network of $m \in \mathbb{N}$ nodes labeled by the set $\mathcal{V}=\lbrace 1,2,...,m \rbrace $ is considered. The topology of the interconnections among nodes is \textit{not} fixed but defined by  a set of graphs $\mathcal{G}(\omega^{*})=(\mathcal{V},\mathcal{E}(\omega^{*}))$ where $\mathcal{E}(\omega^{*})$ is the ordered edge set $\mathcal{E}(\omega^{*}) \subseteq \mathcal{V} \times \mathcal{V}$ and $\omega^{*} \in \Omega^{*}$ where $\Omega^{*}$ is the set of all possible communication graphs, i.e., $\Omega^{*}=\lbrace \mathcal{G}_{1}, \mathcal{G}_{2}, ..., \mathcal{G}_{\bar{N}}  \rbrace$. We assume that $(\Omega^{*},\sigma)$ is a measurable space where $\sigma$ is the $\sigma$-algebra on $\Omega^{*}$. We write $\mathcal{N}_{i}^{in} (\omega^{*})/\mathcal{N}_{i}^{out} (\omega^{*})$ for the labels of agent $i$'s in/out neighbors at graph $\mathcal{G}(\omega^{*})$ so that there is an arc in $\mathcal{G}(\omega^{*})$ from vertex $j/i$ to vertex $i/j$ only if agent $i$ receives/sends information from/to agent $j$. We write $\mathcal{N}_{i}(\omega^{*})$ when $\mathcal{N}_{i}^{in}(\omega^{*})=\mathcal{N}_{i}^{out}(\omega^{*})$. We assume that there are no self-looped arcs in the communication graphs. We also assume that there is no communication delay or noise in delivering a message from agent $j$ to agent $i$.   
	
	Now, we consider the distributed resource allocation problem with $m$ agents that share their local resources. For each node $i \in \mathcal{V}$, there is a private convex cost function $f_{i}:\Re^{n} \longrightarrow \Re$ and a private finite amount of a local resource, given by an $n$-dimensional vector, owned by the agent $i$, i.e., $R_{i} \in \Re^{n}$, that are known to node $i$. The objective of each agent is to collaboratively seek the solution of the following optimization problem using local information exchange with the neighbors and \textit{possibly} switching communication topologies \textit{with} or \textit{without} asynchronous protocol:
	
	\begin{equation}\label{1}
		\begin{aligned}
			& \underset{x_{1}, \hdots, x_{m}}{\text{min}}
			& & \sum_{i=1}^{m} f_{i}(x_{i}) \\
			& \text{subject to}
			& & \sum_{i=1}^{m} x_{i}= \sum_{i=1}^{m} R_{i}
		\end{aligned}
	\end{equation}
	where $x_{i} \in \Re^{n}$ is the decision variable of agent $i$, and the constraint represents resource sharing among the agents.

	So far, we have introduced the distributed \textit{resource allocation} optimization (\ref{1}). Another problem in multi-agent systems is distributed \textit{consensus-based} optimization where the objective of each agent is to collaboratively seek a solution of the following optimization problem:
	\begin{equation}
		\begin{aligned}
			& \underset{s}{\text{min}}
			& & \sum_{i=1}^{m} h_{i}(s) 
		\end{aligned}
	\end{equation}
	where $s \in \Re^{n}$, and $h_{i}:\Re^{n} \longrightarrow \Re$ is the private cost of agent $i$. As a matter of fact, the formulation of distributed consensus-based optimization by using local variables of agents is:
	\begin{equation}\label{distconsensusbased}
		\begin{aligned}
			& \underset{x_{1}, \hdots, x_{m}}{\text{min}}
			& & \sum_{i=1}^{m} h_{i}(x_{i}) \\ 
			& \text{subject to}
			& & x_{1}=x_{2}=\hdots=x_{m}
		\end{aligned}
	\end{equation}
	where $x_{i} \in \Re^{n}$ is the decision variable of agent $i$, and the constraint is achieved by interacting with neighbors and \textit{possibly} switching communication topologies \textit{with} or \textit{without} asynchronous protocol. Relevantly, the set 
	\begin{equation}\label{consensussubspace}
		\mathcal{C}:=\{ x \in \Re^{mn} | x_{i}=x_{j}, 1 \leq i,j \leq m, x_{i} \in \Re^{n} \}
	\end{equation}
	is known as \textit{consensus subspace}. 
	

	In Section IV, we show that \textit{first-order} gradient algorithms of distributed consensus-based optimization can also be used to find an optimal solution of (\ref{1}) over static directed networks with synchronous protocol under \textit{weaker} assumptions (than strong convexity) such as strict convexity, Lipschitz gradient, or coercivity.	In Section V, we show that \textit{first-order} gradient algorithms of distributed resource allocation  can also be utilized to find an optimal solution of (\ref{distconsensusbased}) under certain assumptions. 
	
	\textbf{Remark 1:} All \textit{assumptions} given in this paper can lead to derive distributed algorithms for both consensus-based optimization and resource allocation over \textit{static}, \textit{time-varying}, or \textit{random} \textit{directed} networks \textit{with} or \textit{without} asynchronous protocols under \textit{arbitrary} initialization.
	
	\section{Solving (\ref{1}) by Using Distributed Consensus-Based Optimization Algorithms}

	The Lagrangian $\mathcal{L}:\Re^{mn} \times \Re^{n}\longrightarrow \Re$ associated with (\ref{1}) (see \cite{boydbook}) is defined as

	\begin{equation}\label{lagrangian}
		\mathcal{L}(x,y):= \sum_{i=1}^{m} f_{i}(x_{i}) +\sum_{i=1}^{m} y^{T}(R_{i}-x_{i})
	\end{equation}
	where $x=[x_{1}^{T}, \hdots, x_{m}^{T}]^{T}$, and $y$ is the \textit{Lagrange multiplier} associated with the equality constraint in (\ref{1}). It is known that \textit{Lagrange's dual} problem of (\ref{1}) is distributed consensus-based optimization (e.g., see \cite[Sec. III.A]{Zhangconsensusresource} for proof), i.e., 
	\begin{equation}\label{distributed}
		\begin{aligned}
			& \underset{y_{1}, \hdots, y_{m}}{\text{min}}
			& & G(Y):=\sum_{i=1}^{m} G_{i}(y_{i}) \\
			& \text{subject to}
			& & y_{1}=y_{2}=\hdots=y_{m}
		\end{aligned}
	\end{equation}
	where $Y=[y_{1}^{T}, \hdots, y_{m}^{T}]^{T}$, and
	\begin{equation}\label{giprivate}
		G_{i}(y_{i}):=f^{*}_{i}(y_{i})- y_{i}^{T} R_{i}
	\end{equation}
	is the private cost of agent $i$, in which
	\begin{equation}\label{fidefinition}
		f^{*}_{i}(y_{i}):=\sup_{x_{i}} (y^{T}_{i} x_{i}-f_{i}(x_{i}))
	\end{equation}
	is the Fenchel conjugate of $f_{i}(x_{i})$ (see Definition 2). 
	
	\textbf{Comment 1:} As seen in (\ref{giprivate}) and (\ref{fidefinition}), the main \textit{difficulty} in solving (\ref{distributed}) is that the functions $G_{i}(y_{i}), i=1, \hdots,m,$ are not \textit{explicit} functions of dual variables $y_{i}$. In fact, $G_{i}(y_{i}), i=1, \hdots,m,$ are derived from optimization problems in (\ref{fidefinition}) whose \textit{gradients} are hard to compute directly. This poses a big challenge to solve (\ref{1}) by using distributed consensus-based optimization algorithms where the cost functions must be explicit functions of the variables. As mentioned in Section I, several researchers have addressed this issue by using \textit{first-order} gradient distributed consensus-based optimization algorithms over \textit{static} networks with synchronous protocol under certain assumptions.

	In this section, we propose a \textit{procedure} (i.e., \textbf{P1}-\textbf{P3} below) that can be used for applying first-order gradient algorithms of distributed consensus-based optimization to solve distributed resource allocation optimization (\ref{1}) under some assumptions. The assumptions given here are \textit{weaker} than those of the existing results and thus offer more algorithm options for using existing distributed consensus-based optimization (see Surveys \cite{survey2019}-\cite{nedic}) which work over \textit{static}, \textit{time-varying}, or \textit{random directed} networks \textit{with} or \textit{without} asynchronous protocols (see also Remark 1).

	For using \textit{first-order} gradient algorithms of distributed consensus-based optimization, we are required to know the properties of first-order gradients of $G_{i}(y_{i}), i=1, \hdots, m$. Consequently, we need to know properties of first-order gradient of $f^{*}_{i}(y_{i}).$ \textit{Therefore, it is critical to know relations between properties of $f_{i}(x_{i})$ and $f^{*}_{i}(y_{i})$.} \textit{These relations are given in Lemmas 1 and 2 below. }

	\textbf{Lemma 1} \cite[Cor. 4.1.4, Th. 4.2.1, Th. 4.2.2, Cor. 4.2.10]{dualitybook}:
	\begin{itemize}
		\item  Let $f_{i}:\Re^{n} \longrightarrow \Re, i=1, \hdots, m,$ be strictly convex, differentiable, and $1-$coercive. Then
		
		\textit{(i)} $f^{*}_{i}$ is likewise finite-valued on $\Re^{n},$ strictly convex, differentiable, and $1-$coercive;
		
		\textit{(ii)} The continuous mapping $\nabla f_{i}$ is one-to-one from $\Re^{n}$ onto $\Re^{n},$ and its inverse is continuous;
		
		\textit{(iii)} We have
		\begin{equation} \label{gradientfstar1}
			f^{*}_{i}(y_{i})=\langle y_{i}, (\nabla f_{i})^{-1} (y_{i}) \rangle - f_{i} ((\nabla f_{i})^{-1} (y_{i})),
		\end{equation}  
		$$\forall y_{i} \in \Re^{n}.$$ 
		
		\item  Assume that $f_{i}:\Re^{n} \longrightarrow \Re, i=1, \hdots, m,$ are $\rho-$ strongly convex where $\rho >0$. Then $dom f^{*}_{i}=\Re^{n}$, and $\nabla f^{*}_{i}$ is $\frac{1}{\rho}-$Lipschitz, i.e.,
		$$\Vert \nabla f^{*}_{i}(s_{1})-\nabla f^{*}_{i}(s_{2}) \Vert \leq \frac{1}{\rho} \Vert s_{1}-s_{2} \Vert, \quad{} \forall s_{1},s_{2} \in \Re^{n}. $$
		
		\item  Let $f_{i}:\Re^{n} \longrightarrow \Re, i=1, \hdots, m,$ be convex and have $L-$Lipschitz gradient on $\Re^{n}$ where $L>0$. Then $f^{*}_{i}$ is $\frac{1}{L}-$strongly convex on each subset $C \subset dom \partial f^{*}_{i}.$ 
		
		\item   Let $f_{i}:\Re^{n} \longrightarrow \Re, i=1, \hdots, m,$ be convex, twice differentiable, and $1-$coercive. Assume, moreover, that $\nabla^{2} f_{i}(x_{i})$ is positive definite for all $x_{i} \in \Re^{n}.$ Then $f^{*}_{i}$ enjoys the same properties and
		\begin{equation}\label{hessianfstar}
			\nabla^{2} f^{*}_{i}(y_{i})=[\nabla^{2} f_{i}(\nabla f^{-1}_{i}(y_{i}))]^{-1}, \quad{} \forall y_{i} \in \Re^{n}.
		\end{equation}
	\end{itemize}

	\textbf{Lemma 2} \cite{conjugategradient}: Let $f_{i}:\Re^{n} \longrightarrow \Re, i=1, \hdots, m,$ be strictly convex.Then 
	\begin{equation}\label{gradientfstar2}
		\nabla f^{*}_{i}(y_{i})=\arg\max_{x_{i}} (y^{T}_{i} x_{i}-f_{i}(x_{i})).
	\end{equation}

	\textbf{Remark 2:} Gradient and Hessian of the Fenchel conjugate $f^{*}_{i}(y_{i}), i=1, \hdots,m,$ are given \textit{explicitly} in (\ref{gradientfstar2}) and (\ref{hessianfstar}), respectively, under the assumptions in Lemma 2 and Lemma 1. One can also compute gradient and Hessian of $f^{*}_{i}$ by using (\ref{gradientfstar1}) under the assumptions in Lemma 1.

	\textbf{Observation 1:} By reformulating distributed resource allocation (\ref{1}) as distributed consensus-based optimization (\ref{distributed}) on Lagrange's dual variables $y_{i}, i=1, \hdots,m,$ and using properties of the Fenchel conjugate functions $f^{*}_{i}(y_{i})$ (mentioned in Lemma 1), \textit{any} first-as well as second-order gradient distributed consensus-based optimization algorithms can be used for solving (\ref{1}) under the assumptions in Lemma 1. Here, we give more assumptions than those of the literature such that these assumptions allow us to utilize distributed consensus-based optimization algorithms over \textit{directed} \textit{static}, \textit{time-varying}, or \textit{random} networks \textit{with} or \textit{without} asynchronous protocol under \textit{arbitrary} initialization (see also Remark 1).
	
	It remains to show how to apply existing \textit{first-order} distributed consensus-based optimization algorithms by using Lemma 1 to derive first-order distributed algorithm for resource allocation (\ref{1}). 
	
	\textit{The procedure to solve (\ref{1}) by solving (\ref{distributed}) is given in} \textbf{P1}-\textbf{P3} \textit{as follows:}
	
	\textbf{P1}. Consider distributed consensus-based optimization (\ref{distributed}) with variables $y_{i}$ and select an existing first-order gradient distributed algorithm and its related theorem for convergence. 
	
	\textbf{P2}. By using the selected distributed algorithm and conditions of the related theorem, impose suitable assumptions on first-order gradient of $f_{i}^{*}(y_{i})$ defined in (\ref{fidefinition}) (see Remark 2); also impose suitable assumptions on $f_{i}(x_{i})$ by using Lemmas 1 and 2; consequently, compute the first-order gradient of $G_{i}(y_{i})$ defined in (\ref{giprivate}) by using Lemmas 1 and 2; therefore, the optimization (\ref{distributed}) is solved by the chosen algorithm whose conditions are satisfied, namely all agents get consensus on a \textit{dual} optimal solution $y^{*}$ of (\ref{distributed}). 
	
	\textbf{P3}. Each agent obtains the \textit{primal} optimal solution $x_{i}^{*}, i=1, \hdots,m,$ of (\ref{1}) via $x_{i}^{*}=\displaystyle \arg\max_{x_{i}} (x_{i}^{T} y^{*}-f_{i}(x_{i}))$ if exists. Consequently, the relevant theorem for solving (\ref{1}) can be developed from the related theorem of the selected algorithm in \textbf{P1}.


	In the following subsection, based on Observation 1, we apply an existing \textit{first-order} gradient distributed discrete-time algorithm for solving distributed consensus-based optimization to derive a \textit{first-order} gradient discrete-time algorithm for solving distributed resource allocation (\ref{1}). We explicitly show how to use \textbf{P1}-\textbf{P3} above to derive such algorithms and to develop related theorems. As we mentioned in Contribution in Section I, we only develop \textit{first-order} distributed algorithms in this paper.

	\subsection{A First-Order Totally Asynchronous Algorithm for Solving (\ref{1})}

	To the best of authors' knowledge, the discrete-time algorithm presented in this subsection is the \textit{first} totally asynchronous (see footnote 2) algorithm for solving (\ref{1}).
	
	Let the \textit{weighted matrix of the graph} (see also footnote 3) be defined as $\mathcal{W}(\omega^{*})=[\mathcal{W}_{ij}(\omega^{*})]$  with $\mathcal{W}_{ij}(\omega^{*})=a_{ij}(\omega^{*})$ for $j \in \mathcal{N}_{i}^{in}(\omega^{*}) \cup \{ i \}$, and $\mathcal{W}_{ij}(\omega^{*})=0$ otherwise, where $a_{ij}(\omega^{*})>0$ is the scalar constant weight that agent $i$ assigns to the information $x_{j}$ received from agent $j$. For instance, if $\mathcal{W}(\mathcal{G}_{\tilde{k}})=I_{m}$, for some $1 \leq \tilde{k} \leq \bar{N}$, implies that there are no edges in $\mathcal{G}_{k}$, or/and all nodes are not activated for communication updates in asynchronous protocol or both. The advantage of this formulation is that it \textit{unifies} switching networks with asynchrony in a better understandable formulation (see \cite{alavianiTAC}, \cite{alavianisignalprocessing} for details).

	\textbf{Assumption 1} \cite{alavianiTAC}: 
	\begin{itemize}
		\item \textbf{A1}. The weighted graph matrix $\mathcal{W}(\omega^{*})$ is doubly stochastic for each $\omega^{*} \in \Omega^{*}$ where $\Omega^{*}$ is defined in Section III, i.e.,
		
		\textit{i)} $\sum_{j \in \mathcal{N}_{i}^{in}(\omega^{*}) \cup \{ i \}} \mathcal{W}_{ij}(\omega^{*})=1, i=1,2,...,m$,
		
		\textit{ii)} $\sum_{j \in \mathcal{N}_{i}^{out}(\omega^{*}) \cup \{ i \}} \mathcal{W}_{ij}(\omega^{*})=1, i=1,2,...,m.$
		
		\item \textbf{A2}. The union of all of the graphs in $\Omega^{*}$ is strongly connected, i.e., $Re[\lambda_{2}(\sum_{\omega^{*} \in \Omega^{*}} (I_{m}-\mathcal{W}(\omega^{*})))]>0$.
		
		\item \textbf{A3}. There exists a nonempty subset $\tilde{K} \subseteq \Omega^{*}$ such that $FVP(T)=\{ \tilde{z} | \tilde{z} \in \Re^{mn}, \tilde{z}=T(\bar{\omega},\tilde{z}), \forall \bar{\omega} \in \tilde{K} \}$ where $T(\omega^{*},x):=W(\omega^{*})x$, $\forall \omega^{*} \in \Omega^{*}$, $W(\omega^{*}):=\mathcal{W}(\omega^{*}) \otimes I_{n}$, and each element of $\tilde{K}$ occurs infinitely often almost surely.
	\end{itemize}
	
	\textbf{Remark 3} \cite{alavianiTAC}: If the sequence $\{ \omega^{*}(t) \}_{n=0}^{\infty}$ is mutually independent with $\sum_{t=0}^{\infty} Pr_{t}(\bar{\omega})=\infty$ where $Pr_{t}(\bar{\omega})$ is the probability of $\bar{\omega}$  occurring at time $t$, then part \textbf{A3} in Assumption 1 is satisfied. Moreover, any ergodic stationary sequences $\{ \omega^{*}(t) \}_{t=0}^{\infty}, Pr(\bar{\omega})>0,$ satisfy part \textbf{A3} in Assumption 1. Consequently, any time-invariant Markov chain with its unique stationary distribution as the initial distribution satisfies part \textbf{A3} in Assumption 1. Also B-connectivity (see footnote 1) assumption satisfies part \textbf{A3} in Assumption 1.
	
	\textbf{Assumption 2:} $f_{i}:\Re^{n} \longrightarrow \Re, i=1, \hdots, m,$ are $\frac{1}{K}$-strongly convex where $K>0$, and $\nabla f_{i}(x_{i})$ are $\frac{1}{\mu}$-Lipschitz continuous where $\mu >0$. 
	
	In Theorem 1 below and its proof, we show how to apply \textbf{P1}-\textbf{P3} and Observation 1 to derive the \textit{first} totally asynchronous first-order gradient distributed discrete-time algorithm for resource allocation (\ref{1}).

	\textbf{Theorem 1:} Consider distributed consensus-based optimization (\ref{distributed}) with Assumptions 1 and 2. Let $\beta \in (0,\frac{2 \mu}{K^{2}})$ and $\alpha(t) \in [0,1], t \in \mathbb{N} \cup \lbrace 0 \rbrace,$ satisfy
	
	\textit{(a)} $\displaystyle \lim_{t \longrightarrow \infty} \alpha(t)=0,$
	
	\textit{(b)} $\sum_{t=0}^{\infty} \alpha(t)=\infty.$

	Then starting from any initial points, the sequences $\{y_{i}(t)\}_{t=0}^{\infty}, i=1, \hdots,m,$ generated by     
	\begin{align}
		y_{i}(t+1)&=\alpha(t) (y_{i}(t)- \beta z_{i}(t))+(1-\alpha(t)) ((1-\eta) y_{i}(t) \nonumber \\
		&\quad{} \quad{} +\eta \sum_{j \in \mathcal{N}_{i}^{in}(\omega^{*}(t)) \cup \{ i \}} \mathcal{W}_{ij}(\omega^{*}(t)) y_{j}(t)), \label{resourcealgorithm1}\\
		z_{i}(t)&=x_{i}(t)-R_{i}, \label{resourcealgorithm12}\\
		x_{i}(t)&=\arg\min_{q_{i}} (f_{i}(q_{i})-y^{T}_{i}(t) q_{i}), \label{resourcealgorithm2}
	\end{align}
	where $\eta \in (0,1)$, and $\omega^{*}(t)$ is a realization of $\Omega^{*}$ at time $t$, globally converge almost surely and in mean square to the unique solution of dual problem (\ref{distributed}). Moreover, the sequences $\{x_{i}(t)\}_{t=0}^{\infty}, i=1, \hdots,m,$ converge almost surely and in mean square to the unique primal solution of distributed resource allocation (\ref{1}). 
	
	\textit{Proof:} Let's apply \textbf{P1}-\textbf{P3} and Observation 1 to prove Theorem 1.

	\textbf{P1}: \textit{We consider the first totally asynchronous algorithm proposed in \cite{alavianiTAC} for solving distributed consensus-based optimization (\ref{1}) and its related theorem for convergence (i.e., Theorem 3 in Appendix A).}

	The first totally asynchronous distributed discrete-time algorithm proposed in \cite{alavianiTAC} in compact form is given as follows:
	\begin{align}
		Y(t+1)&=\alpha(t) (Y(t)- \beta \nabla G(Y(t))) \nonumber\\
		&\quad{} +(1-\alpha(t)) ((1-\eta) Y(t)+\eta W(\omega^{*}(t))Y(t)) \label{TACalgorithm}
	\end{align}
	where $W(\omega^{*}):=\mathcal{W}(\omega^{*}) \otimes I_{n}$, $\eta \in (0,1),$ and $\omega^{*}(t)$ is a realization of $\Omega^{*}$ at time $t$. From Theorem 2 in Appendix A, Algorithm (\ref{TACalgorithm}) converges almost surely and in mean square to the unique optimal solution $y^{*}$ of (\ref{distributed}) under Assumption 1 above and Assumption 3 in Appendix A (see Appendix A for details).

	\textbf{P2}: \textit{We impose suitable assumptions on the first-order gradient of $G_{i}(y_{i})$ defined in (\ref{giprivate}) and calculate it by using Lemmas 1-2.}
	
	Since, from Assumption 3 in Appendix A, we require that $G_{i}(y_{i}), i=1, \hdots, m,$ be $\mu$-strongly convex, and $\nabla G_{i}(y_{i})$ be $K$-Lipschitz continuous, we need from Lemma 1 to impose Assumption 2 on $f_{i}(x_{i})$. Therefore, based on Lemma 2 and Assumption 2, we can compute the gradient of $G_{i}(y_{i})$ as follows:
	\begin{eqnarray}\label{gradientGi}
		\nabla G_{i}(y_{i})=\arg\max_{x_{i}} (y^{T}_{i} x_{i}-f_{i}(x_{i})) - R_{i}.
	\end{eqnarray}
	Now to develop Theorem 1 from Theorem 2 in Appendix A, we substitute Assumption 2 for Assumption 3 since Assumption 2 is on primal variables $x_{i}, i=1, \hdots,m.$
	
	\textbf{P3}: \textit{Each agent obtains the primal optimal solution $x_{i}^{*}, i=1, \hdots,m,$ of (\ref{1}) via $x_{i}^{*}=\displaystyle \arg\max_{x_{i}} (x_{i}^{T} y^{*}-f_{i}(x_{i}))$. }

	We obtain the primal optimal solutions $x_{i}^{*}, i=1, \hdots,m,$ as 
	\begin{equation}\label{argminGi}
		x_{i}^{*}=\displaystyle \arg\min_{x_{i}} (f_{i}(x_{i})-x_{i}^{T} y^{*})
	\end{equation}
	that exist from Assumption 2. 
	
	Therefore, we obtain Algorithm (\ref{resourcealgorithm1}) from (\ref{TACalgorithm})-(\ref{argminGi}). Note that Algorithm (\ref{TACalgorithm}) is in a compact form and can be viewed based on local information as Algorithm (\ref{resourcealgorithm1}). As a matter of fact, the distributed algorithm (\ref{resourcealgorithm1}) can converge almost surely and in mean square to the primal optimal solutions $x_{i}^{*}, i=1, \hdots,m,$ of (\ref{1}) under Assumptions 1-2. As seen in Theorem 2 in Appendix A, the convergence occurs over infinite-time horizon. Thus the proof of Theorem 1 is complete.

	\textbf{Remark 4:} Algorithm (\ref{resourcealgorithm1}) does \textit{not} require \textit{a priori} B-connectivity or distribution assumption on switching communication graphs for convergence that results in the \textit{first} totally asynchronous algorithm for solving distributed \textit{resource allocation} (\ref{1}) (see also Remark 3). In general, the rate of convergence cannot be established for a totally asynchronous algorithm. However, determining rate of convergence of Algorithm (\ref{resourcealgorithm1}) under suitable assumptions remains a part of the future work.

	\textbf{Remark 5:} We refer the interested reader to \cite{alavianisecondoreder} for deriving the \textit{first} distributed second-order continuous-time algorithm for resource allocation.

	\section{Solving (\ref{distconsensusbased}) by Using Distributed Resource Allocation Algorithms}

	\textit{Fenchel's dual} problem \cite{dualitybook} of distributed consensus-based optimization (\ref{distconsensusbased}) is the following resource allocation (see \cite[Sec. III.A]{Fencheldualityproof} for proof): 
	
	\begin{equation}\label{11}
		\begin{aligned}
			& \underset{y_{1}, \hdots, y_{m}}{\text{min}}
			& & \sum_{i=1}^{m} h_{i}^{*}(y_{i}) \\
			& \text{subject to}
			& & \sum_{i=1}^{m} y_{i}= \textbf{0}_{n}
		\end{aligned}
	\end{equation}
	where $h_i^{*}(y_{i})$ is Fenchel conjugate function of $h_{i}(x_{i})$ (see Definition 2). We require to mention that similar challenge to the challenge explained in Comment 1 holds for (\ref{11}). Similar to Observation 1, we have the following observation for solving (\ref{distconsensusbased}) by using distributed resource allocation's algorithms.
	
	\textbf{Observation 2:} By reformulating distributed consensus-based optimization (\ref{distconsensusbased}) as distributed resource allocation (\ref{11}) on \textit{Fenchel's dual} variables $y_{i}, i=1, \hdots,m,$ and using properties of the Fenchel conjugate functions $h^{*}_{i}(y_{i})$ presented in Lemmas 1-2, \textit{any} first-as well as second-order gradient distributed resource allocation algorithms can be used for solving (\ref{distconsensusbased}) under the assumptions of Lemma 1. Indeed, Remark 1 is also valid here.

	
	As stated previously, the process for deriving first-order gradient distributed consensus-based optimization algorithms from distributed resource allocation algorithms is similar to the procedures of the proof of Theorems 1 in the previous section.

	\section{Conclusions}
	In this paper, we have proposed a procedure to utilize distributed consensus-based optimization algorithms for solving resource allocation optimization and vice versa over switching networks with/without synchronous protocol. We have shown that first-order gradient algorithms of distributed consensus-based optimization  can be used to find an optimal solution of distributed resource allocation under weaker assumptions than those of the literature over static networks under synchronous protocol. This offers more choices for distributed consensus-based optimization algorithms to be utilized for solving distributed resource allocation problem that can overcome difficulties of existing results. We have also shown that distributed resource allocation algorithms can be employed to derive distributed consensus-based optimization algorithms. This paper has presented the first totally asynchronous algorithm for solving distributed resource allocation by using the proposed procedure. These results can also be utilized to derive distributed algorithms for static, time-varying, or random directed networks with or without synchronous protocol with arbitrary initialization.

	\appendices
	\section{}
	
	The following theorem has been given in \cite{alavianiTAC} for distributed optimization (\ref{distributed}).

	\textbf{Assumption 3} \cite{alavianiTAC}: $G_{i}(y_{i}), i=1, \hdots, m,$ is $\mu$-strongly convex, and $\nabla G_{i}(y_{i})$ is $K$-Lipschitz continuous.

	\textbf{Theorem 2} \cite[Cor. 1-2]{alavianiTAC}: Consider dual problem (\ref{distributed}) with Assumptions 1 and 3. Let $\beta \in (0,\frac{2 \mu}{K^{2}})$ and $\alpha_{n} \in [0,1], n \in \mathbb{N} \cup \lbrace 0 \rbrace$, satisfy \textit{(a)} and \textit{(b)} in Theorem 1. Then starting from any initial point, the sequence generated by the distributed algorithm (\ref{TACalgorithm}) globally converges almost surely and in mean square to the unique solution of (\ref{distributed}).

	\textbf{Remark 6:} An example of the diminishing step size $\alpha_{n}$ satisfying \textit{(a)} and \textit{(b)} in Theorem 1 is $\alpha(t) :=\frac{1}{(1+t)^{\zeta}}$ where $\zeta \in (0,1]$.

	\textbf{Remark 7:} As seen from part \textbf{A3} in Assumption 1, Algorithm (\ref{TACalgorithm}) does \textit{not} require a priori B-connectivity or distribution assumption on switching communication graphs for convergence that results in the \textit{first} totally asynchronous algorithm for solving distributed consensus-based optimization (\ref{distributed}) (see \cite{alavianiTAC} for details).

	

	%
	

\end{document}